\documentclass{article}
\usepackage{graphicx}
\usepackage{amsmath}
\usepackage{amssymb}
\usepackage{hhline}
\usepackage{epsfig}

\renewcommand\a{\alpha}
\renewcommand\b{\beta}
\newcommand\bff[1]{\underline{\bf #1}$^\circ$)}
\newcommand\B{{\cal B}}
\newcommand\BB{{\hat\B}}

\renewcommand\c{\circ}
\newcommand\C{{\mathbb C}}
\newcommand\CC{{\cal C}}
\newcommand\cyl{{\cal C}}
\newcommand\cli{\hat{\cal C}}
\renewcommand\d{\delta}

\newcommand\DD{{\cal D}}
\newcommand\e{\eta}

\newcommand\g{\gamma}

\newcommand\I{{\cal I}}
\newcommand\J{{\hat J}}

\renewcommand\l{\lambda}

\newcommand\oo{\omega}
\newcommand\ot{\otimes}
\newcommand\ok{{\cal O}}

\newcommand\op[1]{\mathop{\rm #1}\nolimits}
\newcommand\p{\partial}
\newcommand\po{$\!\!\!{\bf .}$ }
\renewcommand\P{\Phi}
\newcommand\PP{{\cal P}}

\newcommand\R{{\mathbb R}}

\renewcommand\t{\times}
\newcommand\te{\theta}
\newcommand\ti{\tilde}

\newcommand\ve{\varepsilon}
\newcommand\vp{\varphi}

\newcommand\vph{\vphantom{$\frac{2^2}{2_2}$}}

\newcommand\we{\wedge}
\newcommand\x{\xi}

\newcommand\z{\sigma}
\newcommand\Z{{\mathbb Z}}

\newcounter{acb}

\def\Rom#1{\uppercase\expandafter{\romannumeral#1}}

\newcommand\1{{\bf 1}}

\newcommand\qed{\phantom{\underline{y}}\hfill\hfill$\square$}

\newcommand\bib[1]{\bibitem[#1]{#1}}

 \newcommand{\ddfrac}[2]{\dfrac{\mathstrut\displaystyle #1}
 {\mathstrut\displaystyle #2}}

\newtheorem{theorem}{Theorem}

\newtheorem{theo}{Theorem}
\newtheorem{theu}{Theorem}

\newtheorem{cor}[theorem]{Corollary}
\newtheorem{dfn}[theo]{Definition}
\newtheorem{lem}[theorem]{Lemma}
\newtheorem{prop}[theorem]{Proposition}
\newtheorem{rk}[theu]{Remark}
\newenvironment{proof}[1][Proof]{\textbf{#1.} }{\qed}

\newenvironment{ex}{\trivlist \item[\hskip \labelsep{\bf Example.}]}
{\endtrivlist}

\newcounter{f}
\makeatletter
\newcommand{\@thefnmark}{$^\fnsymbol{f}$}
\renewcommand{\@makefnmark}{\hbox{\mathsurround=0pt
                           $^{\fnsymbol{f}}$}}
\renewcommand{\@makefntext}[1]{\parindent=1em\noindent
            \hbox to 1.8em{\hss$^{\fnsymbol{f}}$}#1}
\makeatother

\begin{document}

\title{{\bf Nijenhuis tensors in \\
pseudoholomorphic curves neighborhoods}}
\author{Boris~S.\ Kruglikov}
 \date{December 2000}
\maketitle

 \begin{abstract}
In this paper%
 \footnote{Author's research is supported by grants INTAS 96-0713,
NWO-RFBR 047-008-005.}
 we consider the normal forms of almost complex structures
in a neighborhood of pseudoholomorphic curve.
We define normal bundles of such curves and study the properties of
linear bundle almost complex structures.
We describe 1-jet of the almost complex structure along a curve in terms
of its Nijenhuis tensor. For pseudoholomorphic tori we investigate the
problem of pseudoholomorphic foliation of the neighborhood. We get some
results on nonexistence of the tori deformation.
 \end{abstract}


\section*{Introduction}

\hspace{13.5pt}
Let $M^{2n}$ be an almost complex manifold, i.e.\ it is equipped with the
tensor field $J\in T^*M\ot TM$ such that $J^2=-\1$. A submanifold
$N\subset M$ is called {\it pseudoholomorphic\/} (PH-submanifold) if its
tangent bundle $TN\subset TM$ is invariant under the operator $J$. In this
paper we study pseudoholomorphic curves, which are 2-dimensional
submanifolds. Generically there are no other PH-submanifolds.

The investigation of pseudoholomorphic curves was initiated by Gromov
in~\cite{G} (but first they appeared in~\cite{NW}).
The structure of the moduli space of such curves plays an
important role in symplectic geometry (\cite{MS}). Generically
nonexceptional PH-spheres occur not discretely but in families. The same
situation is also with PH-disks (\cite{K3}). However it is known
that even in the complex situation generically exist discrete
holomorphic tori (\cite{A1}).

In~\cite{Mo} Moser posed a question what type of KAM-theory can be
constructed for a PH-foliation of an almost complex torus $T^{2n}$. If
this is a foliation by PH-tori $T^2$ then such a theory
can be expected only under very special conditions. Namely
in general position PH-tori are discrete and their number was
investigated by Kuksin in~\cite{Ku1}. Moser considered a
foliation by entire PH-lines $\C\to T^{2n}$ with the slope
of general position. The main result states that under a small
almost complex perturbations of the standard complex structure $J_0$ on
$T^{2n}$ many leaves persist. If the perturbation is big but
tame-restricted then only some of the leaves persist. This was proved by
Bangert in~\cite{B}. Slight generalization of his result with another
approach see in~\cite{KO}.

In section~1 we study the structure of the neighborhood
of a PH-torus in a 4-dimensional manifold. The dimension four is
quite specific for the Nijenhuis tensor and generically it leads
to some canonical $\{e\}${\em-structure\/} (theorem~\ref{th:1!}).
This gives sufficient number of moduli for classification of almost
complex structures in a germ of a point as well as of a PH-curve.

The moduli are absent in the complex situation and
many results of the complex theory fail in the almost complex case.
In particular this is so with Arnold theory of elliptic curves
neighborhoods (which he called Floquet-type theory~\cite{A1}).
In this paper we recover some traces of normal bundles to holomorphic
curves theory in almost complex situation.

The Nijenhuis tensor is determined by 1-jet of the almost complex
structure. 1-jet of a complex structure on a holomorphic curve determines
its holomorphic normal bundle. So to define pseudoholomorphic normal
bundle we should involve the Nijenhuis tensor of the structure $J$.

In section 2 we define abstract linear bundle almost complex structures
and prove that the normal bundle $N_\CC M=TM/T\CC$ to a PH-curve
$\CC$ has a canonical structure of this kind.
Then we show that the Nijenhuis tensor of such structure is
quite special in a sense that its characteristic 2-dimensional
distribution is integrable and the leaves coincide with the fibers.
This allows to find normal forms of linear bundle almost complex
structures and thus normal forms of 1-jet of an almost complex structure
along a PH-curve. Of course the formula involves the Nijenhuis tensor
only because it is a complete invariant of the almost complex structure
(\cite{K1}).

In section 3 we investigate the germ of neighborhoods of a PH-torus
$T^2\subset M^4$ by means of pseudoholomorphic curves.
To get some analogy with the complex situation
we note that for a neighborhood of an elliptic curve
there is a special foliation by cylinders.
Part of these cylinders persists in the almost complex category.
We construct families of cylinders and investigate some
natural questions of transports and monodromy of the transversal
PH-disks.

The last part of section 3 is devoted to the problem of
pseudoholomorphic tori deformation.
It mainly concerns the linear almost complex bundles since
these model normal bundles of pseudoholomorphic curves. We deduce the
equation on the deformed tori, study it and describe some cases when there
exists a unique solution. This means that the torus is isolated and
persistent. The first example is a sort of geometric version of the
Moser's non-deformation example \cite{Mo}.

Since the most of the paper is about normal bundles in almost complex
category we sketch at the appendix what happens in other geometries,
namely Riemannian.

The paper was started by a question of Arnold (\cite{A2}; 1993-25)
about almost complex version of his theory. As we show the answers
to most of his questions are negative (in a sense: "linearizable"
along a curve complex structures form a typical set, while in the
almost complex case, codimension of such "linearizable"
structures is $\infty$. However we describe the obstructions for
"linearizations" -- equivalence to the normal bundle structure -- in
\S\ref{S2.2}). But this does not bring the issue
to a close. For example there are still questions about
pseudoholomorphic {\em foliations\/} of a PH-torus neighborhood by
cylinders (see section 3) to be investigated.
I expect this can be done by a method similar to Moser's \cite{Mo}.

\section{Moduli of germs of pseudoholomorphic curves neighborhoods}

\subsection{Nijenhuis tensor characteristic distribution}

\hspace{13.5pt}
Here we recall basic facts about Nijenhuis tensors in dimension 4.
They generalize theorem 7~\cite{K1}.

Nijenhuis tensor of almost complex manifold $(M,J)$
is the $(2,1)$-tensor $N_J\in\Lambda^2T^*M\ot TM$ given by
the formula:
 \begin{equation}
N_J(\x,\e)=[J\x,J\e]-J[J\x,\e]-J[\x,J\e]-[\x,\e].
 \label{e0}
 \end{equation}

This tensor satisfies $J$-linearity
$N_J(J\x,\e)=N_J(\x,J\e)=-JN_J(\x,\e)$. So it can
be considered as an antilinear map $N_J:\Lambda^2\C^2\to\C^2$,
$\C^2=(T_x M^4,J)$. Thus the image is invariant under $J$
and hence is a complex line $\C\subset\C^2$ (if $N_J\ne0$) or a point $0$
(if $N_J=0$).

Suppose $N_J\ne0$ in a neighborhood of $x\in M$. We consider the above
complex line as $J$-invariant 2-dimensional real subspace of $T_xM$ and
denote it by $\Pi_x$. Thus we obtain a distribution of 2-dimensional
planes in a neighborhood of $x$.

 \begin{dfn}\po
Let us call $\Pi^2=\op{Im}N_J$ the
{\em Nijenhuis tensor characteristic distribution\/} on 4-dimensional
almost complex manifold $(M^4,J)$.
 \end{dfn}

 \begin{rk}\po
At the points, where $N_J=0$, we have $\op{Im}N_J=0$. Therefore we should
call $\Pi$ differential systems. However singular points of $\Pi$, where
the rank is 0, are generically isolated and so are not relevant to our
future discussion. Thus we keep the term distribution.
 \end{rk}

This distribution $\Pi^2$ is in general situation nonintegrable.
Therefore it has nontrivial derivative $\Pi^3=\p\Pi^2$,
that is the differential system with the $C^\infty(M)$-module of sections
$\PP_3=C^\infty(\Pi^3)$ generated by the module of
sections $\PP_2=C^\infty(\Pi^2)$ and its self-commutator:
$\PP_3=[\PP_2,\PP_2]$. $\Pi^3$ is not distribution everywhere,
it can have singularities. These singularities form a (stratified)
manifold $\Sigma_2'$ of $\op{codim}=2$ with singularities.

Consider a point outside $\Sigma'_2$. Generically the distribution $\Pi^3$
is not integrable and so $\PP_4=[\PP_3,\PP_3]=C^\infty(TM)=\DD(M)$
is the module of all vector fields,
so that the next distribution is the whole tangent space. Moreover we can
assume $[\PP_2,\PP_3]=\DD(M)$. This again fails to be so outside some
other (stratified) manifold with singularities $\Sigma_2''$ of
$\op{codim}=2$.

Let $x\notin\Sigma_2'$. Then $\Pi^2_x\subset\Pi^3_x$ has a transversal
measure. Actually there exist vectors $\x_1,\x_2\in\Pi^2_x$,
$\x_3\in\Pi^3_x\setminus\Pi^2_x$ such that $N_J(\x_1,\x_3)=\x_1$,
$N_J(\x_2,\x_3)=-\x_2$ (because a $J$-antilinear isomorphism
$N_J(\cdot,\x_3):\Pi^2_x\to\Pi^2_x$ is orientation reversing).
The fields $\x_1,\x_2$ are defined up to
multiplication by a constant, while $\x_3(\!\!\!\!\mod\Pi^2_x)$ is defined
up to multiplication by $\pm1$. Therefore $\Pi^3/\Pi^2$ is normed.

By a similar reason $T_xM/\Pi^3_x$ is normed outside
$\Sigma_2'\cup\Sigma_2''$.  The transversal measure for $\Pi_x^3\in T_xM$
is given by the vector $\x_4=J\x_3$.

Note that $\Pi^3_x/\Pi^2_x$ is oriented. Actually
$[\x_1,\x_2]\op{mod}\Pi^2_x$
depends only on the values of $\x_1,\x_2$ at the point $x$. It is a
vector $f\x_3\op{mod}\Pi^2_x$ for some $f$.
So if we require $\x_2=J\x_1$ then $\x_3$ can be chosen so that $f>0$.
This produces coorientation of $\Pi^2_x\subset\Pi^3_x$ and then
via $J$ coorientation of $\Pi^3_x\subset T_xM$.

Moreover the requirement $f=1$ determines canonically vector field
$\x_1$ (still however up to $\pm1$) and hence $\x_2=J\x_1$.
Then we set $\x_3=[\x_1,\x_2]$ and $\x_4=J\x_3$.
So the pair $(\x_1,\x_2)$ is defined canonically up to a sign and the pair
$(\x_3,\x_4)$ is absolutely canonical.

 \begin{theorem}\po\label{th:1!}
Let almost complex structure $J$ be of general position. Then
at generic points $x\in M^4$ a canonical frame
$(\x_1,\x_2,\x_3,\x_4)$ is defined. It restores uniquely
the almost complex operator $J$ and the tensor $N_J$ by the tables:

  \begin{center}
{\renewcommand{\arraystretch}{0}%
\begin{tabular}{|c||c|}
\hhline{-||-}
\strut
$X$ & $JX$ \\
\hhline{:=::=:}
\strut $\x_1$ & $\x_2$ \\
\hhline{|-||-|}
\strut $\x_2$ & $-\x_1$ \\
\hhline{|-||-|}
\strut $\x_3$ & $\x_4$ \\
\hhline{|-||-|}
\strut $\x_4$ & $-\x_3$ \\
\hhline{-||-}
\end{tabular}}
\hspace{30pt}
{\renewcommand{\arraystretch}{0}%
\begin{tabular}{|c||c|c|c|c|}
\hhline{-||----}
\strut $N_J({\scriptstyle\uparrow},{\scriptstyle\leftarrow})$
 & $\x_1$ & $\x_2$ & $\x_3$ & $\x_4$ \\
\hhline{:=::====:}
\strut $\x_1$\ & $0$      & $0$     & $\x_1$   & $-\x_2$  \\
\hhline{|-||----|}
\strut $\x_2$\ & $0$      & $0$     & $-\x_2$   & $-\x_1$ \\
\hhline{|-||----|}
\strut $\x_3$\ & $-\x_1$    & $\x_2$   & $0$      & $0$    \\
\hhline{|-||----|}
\strut $\x_4$\ & $\x_2$   & $\x_1$   & $0$      & $0$      \\
\hhline{-||----}
\end{tabular}}
 \end{center}
 \end{theorem}

Note that reducing a structure to the frame (or $\{e\}$-structure)
solves completely the equivalence problem. We recall briefly
the main idea. Consider the moduli of the problem, i.e.\ functions
$c^i_{jk}$ given by the formula $[\x_j,\x_k]=\sum c^i_{jk}\x_i$.
Denote by ${\mathbb A}=\{c^i_{jk}\}$ the space of all invariants and by
$\Phi:M\to {\mathbb A}$ the "momentum map" $x\mapsto\{c^i_{jk}(x)\}$.
Then two equivalent structures
have the same images and this gives an equivalence of the structures in
most cases. For more details see \cite{S}.

 \begin{rk}\po
If both derivatives of $\Pi^2$ are nontrivial, i.e.\
$x\notin\Sigma_2'\cup\Sigma_2''$ then the distribution $\Pi^2$ has a
canonical normal form. It is called {\em Engel distribution\/}
(\cite{BCG}). It is shown in~\cite{K2} that this and other cases are
realized as Nijenhuis tensors characteristic distributions.
 \end{rk}

\subsection{Invariants of a PH-curve neighborhood}

\hspace{13.5pt}
Let $\CC^2$ be a PH-curve. At every point $x\in\CC$ we have two
$J$-invariant planes $T_x\CC^2$ and $\Pi^2_x$. Obviously in general
position they intersects by zero. When these distributions are considered
along the torus they are generically transversal everywhere except a
finite number of points. Denote this set by $\tilde\Sigma_0$. In general
(transversal) case there are also finite sets
$\Sigma_0'=\Sigma_2'\cap\CC$ and $\Sigma_0''=\Sigma_2''\cap\CC$.
The arrangement of these points
 $$
\Sigma_0=\tilde\Sigma_0\cup\Sigma_0'\cup\Sigma_0''\subset\CC
 $$
gives a (finite-dimensional) invariant.

For points $x\in\CC\setminus\Sigma_0'$ we define field of directions
$L^1=T\CC\cap\Pi^3$. The integral curves of this 1-distribution foliate
the set $\CC\setminus\Sigma_0'$ and in general $\CC$ foliates with only
nondegenerate singular points. Denote the number of elliptic points by
$e(L^1)$ and the number of hyperbolical by $h(L^1)$.
 Note that (topologically stable) points of $\tilde\Sigma_0$
are usually regular points of $L^1$.

 \begin{lem}\po
Under $C^1$-small perturbation of the structure $J$ the foliation
$L^1$ has minimal number of singularities,
$\min\{e(L),h(L)\}=0$, $\max\{e(L),h(L)\}=\dfrac12|\chi(L)|$.
In particular for torus $\CC=T^2$ we get foliation
without singularities.
 \end{lem}

 \begin{proof}
Actually this perturbation is $C^0$-small for $\Pi^2$ and large for
$\Pi^3$ and hence $L^1$. To perform it one collects elliptic-hyperbolic
points at pairs (after a small perturbation one can assume there's nothing
more complicated than separatrix self-connection) and kill them one by one.
The possibility of such a perturbation is an easy calculational
argument similar to theorem 5~\cite{K1}.
 \end{proof}

As follows from the previous subsection the foliation $L^1$ is equipped
with orientation and parallel measure outside $\Sigma_0$.
Thus there exists a canonical vector field $v_1$ along $L^1$.
So the curve $\CC$ has many invariants from dynamical systems point of
view, for example winding number of $v_1$ for every two
independent cycles, linearizations at critical points etc.

But since the foliation is equipped with a coorientation and a transversal
measure we have more. Let $v_2$ be positively cooriented, have transversal
measure 1 and be directed along $JL^1$. Then $(v_1,v_2)$ is a canonical
frame outside $\Sigma_0$. Writing
 $$
[v_1,v_2]=\g_1v_1+\g_2v_2.
 $$
we obtain two invariant (under pseudoholomorphic isomorphisms) functions
$\g_1,\g_2$. We can also construct a function $\g_3$ defined by
$Jv_1=\g_3v_2$. These together with the functions $c^i_{jk}$ from the
previous subsection and others form {\em moduli\/} of the
$\CC$-neighborhoods germ.

Note that all constructed invariants are nontrivial. We illustrate this
for the simplest case of arbitrary winding number on the torus.

 \begin{ex}
Consider the manifold $M=T^2(\vp,\psi)\times\R^2(x,y)$ with almost complex
structure given by
 $$
\begin{array}{ll}
J\p_x=\p_y; &
J\p_\vp=\frac{2+\rho y^2}2\p_\psi+\frac{y^2}2\p_\vp+x\p_x;\\
J\p_y=-\p_x; &
J\p_\psi=-\frac{4+y^4}{4+2\rho y^2}\p_\vp-\frac{y^2}2\p_\psi-
\frac{xy^2}{2+\rho y^2}\p_x-\frac{2x}{2+\rho y^2}\p_y,
\end{array}
 $$
where $\rho\in\R$. The torus $\CC=\{x=y=0\}$ is pseudoholomorphic. We
calculate $N_J(\p_x,\p_\vp)=\x$ for $\x=y(\p_\vp+\rho\p_\psi)-\p_y$.
Thus $\Pi^2=\langle\x,J\x\rangle$ is transversal to $\CC$. Since
$[\x,J\x]_{x=y=0}=-\p_\psi+\rho\p_\vp$ we conclude that the
foliation $\langle L^1\rangle$ is given by the equation
$\{\dot\vp=\rho,\dot\psi=-1\}$, whence the winding number is $-\rho$.
 \end{ex}

\subsection{Holomorphic recollections}\label{HolRec}

\hspace{13.5pt}
Let $\CC$ be a holomorphic curve in some complex surface and let it
be not self-linked, $\CC\cdot\CC=0$.
Consider its co-normal bundle defined as the first term of the
following exact sequence:
 $$
0\to \I/\I^2 \to T^*M\vert_\CC \to T^*\CC \to 0,
 $$
where $\I$ is the ideal of holomorphic functions corresponding to the
curve $\I=\{f:\ok\to\C\,|\,f(\CC)=0\}$ and $\ok$ is the germ of
neighborhoods of $\CC\subset M$. The dual sequence
 \begin{equation}\label{defNB}
0 \to T\CC \to TM\vert_\CC \to N_\CC M \to 0
 \end{equation}
defines the normal bundle.

Let the curve be elliptic $\CC=T^2$. This torus is uniquely characterized
by its periods $(2\pi,\nu)$, where $\op{Im}(\nu)\ne0$. The complex number
$\nu$ is determined almost uniquely in the sense that the corresponding
lattice $\Z^2(2\pi,\nu)$ is unique. The torus is given then as the
quotient of $\C$ by the lattice.

All holomorphic topologically trivial 1-dimensional vector bundles over
this elliptic curve are described as follows.
Consider $\C^2$ with coordinates $(z,\vp)$ and for $\l\in\C\setminus\{0\}$,
$\nu\in\C\setminus\R$ make the identifications
 \begin{equation}
(z,\vp)\simeq(z,\vp+2\pi)\simeq(\l z,\vp+\nu).
 \label{e1}
 \end{equation}
The bundle $E\to T^2$ is now given by $(z,\vp)\mapsto\vp$.

Similarly one classifies $S^2$-bundles for which two distinguished
holomorphic sections "zero" and "infinity" are chosen.

Arnold in \cite{A1}~\S27 considers normal bundle $N_{T^2}M^4$ to an
elliptic curve and prove that if the pair $(\l,\nu)$ is normal nonresonant
($\l^n\ne e^{ik\nu}$ plus some diophantine condition for this pair) then a
small neighborhood of $T^2\subset M$ is biholomorphically equivalent to a
neighborhood of the zero section in $N_{T^2}M$.

Note that the number $\lambda$ is defined by 1-jet of the complex
structure on the torus $T^2$. But in almost complex situation 1-jet is
determined by the field of the Nijenhuis tensors $N_J$ along the torus, so
this $\lambda$ is not defined. We are going to define a family of complex
structures $J_{(\l)}$ which have the same values as $J$ at the points of
the torus $T^2$.

\subsection{Complex structures in a PH-torus neighborhood}\label{S14}

\hspace{13.5pt}
Let $\CC\subset M^4$ be a pseudoholomorphic curve. To define some
local normal coordinates in a neighborhood of $\CC$ we prove

 \begin{prop}\po\label{poravstavat'}
Small neighborhood $\ok$ of a PH-curve $\CC\subset M^4$ can be foliated by
transversal PH-disks $D^2$.
\label{prop:1}
 \end{prop}

This actually follows from Nijenhuis-Wolf theorem~\cite{NW} of existence
of small PH-disk in a given direction (which is normal in our case) and
smooth dependence on this direction. We give another proof based on the
(visibly different) idea we exploit later on.

 \begin{proof}
Let us change the almost complex structure $J$ outside a small
neighborhood $\ok$ of $\CC$ in such a manner that it be integrable.
Moreover it can be done so that some bigger neighborhood is isomorphic to
$\CC\t D^2$ near the boundary. Thus we can glue the neighborhood and
obtain the manifold $M_0=\CC\t S^2$ with the almost complex structure
$J'$ being integrable outside the fixed neighborhood of our PH-curve.
(cf.\ prop.~2\cite{K3}).

Now we introduce a symplectic product-structure $\oo=\oo_1\oplus\oo_2$ on
$(M_0,J')$ and note that for sufficiently small neighborhood $\ok$
symplectic structure $\oo$ tames the structure $J'$.
If $\chi(\CC)<2$ the pseudoholomorphic spheres can lie only in
the multiple of the homology class of the second factor
in $M_0=\CC^2\t S^2$. So we seek for them in the
primitive class $[S^2]$. The Gromov theory (\cite{G}) implies that in the
case of general position almost complex structure $J'$ the manifold
$(M_0,J')$ is foliated by PH-spheres. Taking the intersection of this
foliation with the neighborhood $\ok$ we prove the statement.

If the curve is rational $\CC=S^2$ some more delicate arguments are
required (\cite{M2}). If the almost complex structure $J'$ is
not of general position we apply the compactness theorem from~\cite{MS}
for a generic sequence $J_i\to J'$.
 \end{proof}

Now if $z=x+iy$ is a complex coordinate on some transversal PH-disk
we can define coordinates on all close disks by parallel transport
(along specified rays) using
some connection $\nabla$. If the connection is almost complex
$\nabla J=0$ then on every nearby transversal the coordinate is complex.
However we cannot define the coordinates globally in such a way because of
the holonomy. We consider the special case of elliptic curve $\CC=T^2$ (the
case of rational curve is trivial).

Note that the number $\nu$ is defined for $\CC$. Actually almost complex
structure $J$ on any 2-dimensional manifold is integrable and hence the
torus $T^2$ has periods $(2\pi,\nu)$. In particular we have a
(double-periodic) complex coordinate $\vp=\vp_1+i\vp_2$ on the curve.

We extend the coordinate $\vp$ to the neighborhood $\ok$ by the projection
along the disks, i.e.\ the disks are given by $D_\vp=\{\vp=\op{const}\}$.
We introduce now complex coordinate $z=x+iy$ on every transversal disk
in such a way that they specify the same complex structure on $D_\vp$ as
the restriction of $J$ and our PH-torus is $\{z=0\}$.
This coordinates are multivalued because of the holonomy but we can
choose $z$ so that the gluing rules
(monodromies along the cycles) on the torus are given by~(\ref{e1}) with
any prescribed $\l$. Now assuming the coordinates $(z,\vp)$ are complex
we get the complex structure $J_{(\l)}$ in $\ok$.

The previous discussion and Arnold theorem imply that these structures
form an almost complete family in the following sense

 \begin{prop}\po
Let $\ti J$ be a complex structure in a neighborhood of $T^2\subset M^4$
equal to the almost complex structure $J$ at the points of the torus.
Then if the pair $(\l,\nu)$, determined by $\ti J$ on $T^2$, is normal
nonresonant the germ of $\ti J$ is equivalent to the germs of $J_{(\l)}$.
 \label{prop:2}
 \qed
 \end{prop}

\section{Differential geometry of a pseudoholomorphic curve neighborhood}

\subsection{Normal bundle of a pseudoholomorphic curve}\label{AC on NB}

\hspace{13.5pt}
In this section we do not require $\dim M=4$.

In almost complex case the bundle $TM\vert_\CC$ for a PH-curve
$\CC\subset M$ is no longer holomorphic. So one needs to consider
almost complex bundles. A bundle $\pi:(E,J)\to(\CC,J_0)$ is called
almost complex if
 $$
\pi_*(J\x)=J_0(\pi_*\x).
 $$

 \begin{prop}\po
Consider the normal bundle to a PH-curve $\CC\subset M$ given as vector
space by the sequence (\ref{defNB}). There is a canonical almost complex
structure $\J$ on $N_\CC M$ such that $(N_\CC M,\J)$ is an almost
complex bundle.
 \end{prop}

 \begin{proof}
Recall that on almost complex manifold $(M,J)$ there is always an almost
complex connection $\nabla J=0$. Actually if $\nabla'$ is any
linear connection then $\nabla=\frac12\bigl(\nabla'-J\nabla'J\bigr)$ is
also a linear connection, which preserves the structure $J$.
Moreover we can assume that $\nabla$ is minimal, i.e.\ the torsion
of $\nabla$ equals to its antilinear by each argument part
$T_\nabla=T_\nabla^{--}=\frac14N_J$ (\cite{L}).

 \begin{lem}\po
There exists a minimal almost complex connection such that
the curve $\CC$ is a totally geodesic submanifold, i.e.\
parallel transport of any vector $v\in T\CC$ along a path in $\CC$
belongs again to $T\CC$.
 \end{lem}

 \begin{proof}
To see this let us make a gauge transformation $\tilde\nabla=\nabla+A$,
where $A\in\Omega^1(M,\op{end}_\C\tau_M)$ is a 1-form with values in
complex endomorphisms of the tangent bundle. If we require that this
1-form is symmetric $A\in S^2T^*M\ot_\C TM$, then the new connection
$\tilde\nabla$ is also almost complex and minimal.

Let $\x$ be a vector field on $\CC$ with nondegenerate critical points,
$\x\in\DD(\CC)$. Let $\nabla_\x\x=\e\in\DD(M)$.
We define $A(\x,\x)=-\e$,
$A(J\x,\x)=A(\x,J\x)=-J\e$, $A(J\x,J\x)=\e$ and for other vectors somehow
preserving symmetry and $J$-linearity
(near the critical points there's a more work, we need to choose
$\x$ so that $\e$ has a zero of the second order at these critical
points).

Then $\tilde\nabla_\x\x=0$,
$\tilde\nabla_\x J\x=0$. Therefore by minimality
$\tilde\nabla_{J\x}\x=\tilde\nabla_\x
J\x+[J\x,\x]+\frac14N_J(J\x,\x)=[J\x,\x]\in\DD(\CC)$ and also
$\tilde\nabla_{J\x}J\x\in\DD(\CC)$. So $\tilde\nabla$ preserves $T\CC$.
 \end{proof}

Another way to prove this is to introduce trivial connection in
$\ok(\CC)\simeq \CC\times D$ and then to check that procedures of making
connection almost complex and then minimal do not destroy the property of
$\CC$ being totally geodesic. Let us denote this new connection again by
the symbol $\nabla$.

We introduce a connection $\hat\nabla$ to the bundle $N_\CC M$
by means of parallel transport as follows. Let $v=[\te]\in (N_\CC M)_x$
be the class of $\te\in T_xM$ and let $\g(t)\subset\CC$ be a curve,
$\g(0)=x$. Calculate parallel transport $\te(t)$ of $\te$ along $\g(t)$.
We define $v(t)=[\te(t)]$ to be the parallel transport of $v$ along
$\g(t)$. Since $\CC$ is totally geodesic, the definition is correct
($\hat\nabla$-parallel transport of 0 is 0).
Moreover the connection $\hat\nabla$ is $\R$-linear.
So as usual in the theory of generalized connections we conclude that
$\hat\nabla$ is a linear connection.

Let $a\in N_\CC M$ be a point in the normal bundle with projection
$\pi(a)=x$. Let $T_a(N_\CC M)=H_a\oplus V_a$ be the decomposition
by horizontal and vertical subspaces induced by $\hat\nabla$.
The first space $H_a\stackrel{\pi_*}\simeq T_x\CC$ has a canonical complex
structure $J_1$ induced from $J$ by $\pi_*$, while the second
$V_a\simeq T_xM/T_x\CC$ inherits canonical complex structure $J_2$ from
$J$. So we introduce the structure on $N_\CC M$ by the rule
$\hat J=J_1\oplus J_2$.

Let us show that the constructed almost complex structure $\hat J$
does not depend on the choice of minimal connection $\nabla$
preserving $TC$. Let us change the connection $\tilde\nabla=\nabla+A$,
$A\in\Omega^1(M;\op{end}_\C\tau_M)\cap S^2\tau^*_M\ot\tau_M$.
This affects in the change of decomposition
$T_a(N_\CC M)={\tilde H}_a\oplus V_a$, where
${\tilde H}_a=\op{graph}\{\l_a: H_a\to V_a\}$ and $\l_a=A(\cdot)a$.
We state that $\l_a$ is a complex linear map. Actually
 $$
\l_a(Jw)=A(Jw)a=A(a)Jw=JA(a)w=JA(w)a=J\l_a(w).
 $$
So the complex structure $\hat J$ on $T_a(N_\CC M)$ is canonically defined.
 \end{proof}

 \begin{rk}\po
Let $\phi_t:(\CC,J_0)\to(M^4,J)$ be a family of $J$-holomorphic curves of
the same holomorhpic type with $\phi_0=\op{id}$. Then
$\left.\frac{d}{dt}\right|_{t=0}\phi_t$ is the pseudoholomorphic curve in
$N_\CC M$. Moreover we have a one-parametric (scaled) family of PH-curves.
Since $(N_\CC M)_x=T_xM/T_x\CC$ this family does not depend on the
reparametrization $\ti\phi_t=\phi_t\circ g_t$ for
$g_t\in\op{Aut}(\CC,J_0)$ (note that a curve in the bundle bijectively
projected to the base has a natural parametrization). We will use it in
\S\ref{3.4}.
 \end{rk}

 \begin{rk}\po
There is another approach to get an almost complex structure on $N_\CC M$
(\cite{M2}). Consider a foliation of a neighborhood of $\CC\subset M$ to
transversal PH-disks (proposition \ref{poravstavat'}). Let $A_t$ be the
dilation along the disks. We define $J_t=A_tJA_t^{-1}$ and
$\hat J=\lim_{t\to\infty}J_t$. Now the bundle $(N_\CC M,\hat J)$ is almost
complex. However it is difficult to see canonicity and the properties of
this almost complex structure.
 \end{rk}

\subsection{Almost complex bundles}\label{S2.2}

\hspace{13.5pt}
Now we return to the case $\dim M=4$.

 \begin{prop}\po\label{lem:3}
Let $\pi:(E^4,J)\to (\CC^2,J_0)$ be an almost complex bundle over a curve.
Then the characteristic distribution $\Pi^2=\op{Im}N_J$ is integrable and
is tangent to the fibers $F_x=\pi^{-1}(x)$.
 \end{prop}

 \begin{proof}
Actually:
 $$
 \begin{array}{rcl}
\pi_*N_J(\x,\e)
&\!\!\!=\!\!\!&
\pi_*[J\x,J\e]-\pi_*J[\x,J\e]-\pi_*J[J\x,\e]-\pi_*[\x,\e]\\
&\!\!\!=\!\!\!&
[J_0\pi_*\x,J_0\pi_*\e]-J_0[\pi_*\x,J_0\pi_*\e]
-J_0[J_0\pi_*\x,\pi_*\e]-[\pi_*\x,\pi_*\e]\\
&\!\!\!=\!\!\!&
N_{J_0}(\pi_*\x,\pi_*\e)=0.
 \end{array}
 $$

\vspace{-30pt}
 \end{proof}

 \begin{cor}\po
Codimension of the set of almost complex structures, the germs of which on
the PH-curve $\CC\subset M$ are isomorphic to these of the normal bundle
$\CC\subset N_\CC M$, in the set of all almost complex structures is
infinity.
 \qed
 \end{cor}

 \begin{proof}
Actually if the distribution $\Pi^2$ is nonintegrable in a neighborhood of
the PH-curve $\CC$, then a neighborhood $(\ok,J)$ of the curve is not
isomorphic to a neighborhood of the zero section in the normal bundle
$(N_\CC M,\hat J)$.
 \end{proof}

This property is contrary to its analog in the complex category,
see Arnold theorem~\cite{A1} about neighborhoods of elliptic curves
(\S\ref{HolRec}). So extending the category the discussed property becomes
exceptional.

As the following example shows the integrability of $\Pi^2$ is
necessary but by no means sufficient condition on the structure $J$
to be locally isomorphic to its representative on the normal bundle.

 \begin{ex}
Consider the almost complex structure given
on a $T^2(\vp)\t D^2(z)$ with
$\vp=\vp_1+i\vp_2$, $z=x+iy$ by the formula:
 \begin{equation}\label{251000}
J\p_x=\p_y,\
J\p_{\vp_1}=\p_{\vp_2}+A_1\p_{\vp_1}+A_2\p_{\vp_2}+B_1\p_x+B_2\p_y,
 \end{equation}
with $\left.A_i\right|_{T^2}=\left.B_i\right|_{T^2}=0$.
The condition $\Pi^2=T_*\{\vp=\op{const}\}$ is equivalent to the
following PDE system
 \begin{equation}
 \left\{
 \begin{array}{lll}
\ddfrac{\p A_1}{\p y}=A_1\ddfrac{\p A_1}{\p x}-\ddfrac{1+A_1^2}{1+A_2}
\ddfrac{\p A_2}{\p x}&&\\
\ddfrac{\p A_2}{\p y}=(1+A_2)\ddfrac{\p A_1}{\p x}-A_1\ddfrac{\p A_2}{\p x}
&&
 \end{array}
 \right.
  \label{e:ADD}
 \end{equation}

If the functions $A_i$ satisfy this (Cauchy-Kovalevsky) system the
projection along the leaves is given by the formula $(z,\vp)\mapsto\vp$.
Therefore our structure $J$ is projectable iff $A_i\equiv0$. But there are
nonzero solutions of (\ref{e:ADD}). For example:
 $$
A_1=-\dfrac{x}{1+y},\
A_2=-\dfrac{y}{1+y}.
 $$
 \end{ex}

So the space of projectable structures $J$ are of $\op{codim}=\infty$
among the structures with $\Pi^2$ integrable, which are of
$\op{codim}=\infty$ among all almost complex structures with a fixed 0-jet
on the torus $T^2$.

\subsection{Nijenhuis tensor of normal bundles}

\hspace{13.5pt}
Let $N_\CC M$ be a normal bundle of pseudoholomorphic curve
$\CC\subset M^4$. Then at the points $x\in\CC$ two different Nijenhuis
tensors $N_J$ of ambient structure and $N_{\hat J}$ of the structure of
the normal bundle are defined. Are there some relations between these two
tensors? Or between characteristic distributions of these structures?
As the following examples show the answer is negative.

 \begin{ex}
\underline{"Parallel distribution $\Pi^2$"}.
Let $M=\R^4(x_1,y_1,x_2,y_2)$ be equipped with almost complex structure
 $$
J\p_{x_1}=\p_{y_1},\
J\p_{y_1}=-\p_{x_1},\
J\p_{x_2}=\p_{y_2}+x_1\p_{x_1},\
J\p_{y_2}=-\p_{x_2}-x_1\p_{y_1}.
 $$
Then the curve $\CC=\{x_2=y_2=0\}$ is pseudoholomorphic. We calculate
that the following are the Nijenhuis tensor of the structure $J$ and a
minimal almost complex connection
(note that $\Pi^2=\op{Im}N_J=\langle\p_{x_1},\p_{y_1}\rangle$):

 \vspace{1pt}
  \begin{center}
{\renewcommand{\arraystretch}{0}%
  \begin{tabular}{|c||c|c|c|c|}
 \hhline{-||----}
\vph
 \hspace{-5pt}$N_J({\scriptstyle\uparrow},
 {\scriptstyle\leftarrow})$\hspace{-5pt}
 & $\p_{x_1}$ & $\p_{y_1}$ & $\p_{x_2}$ & $\p_{y_2}$ \\
 \hhline{:=::====:}
\vph $\p_{x_1}$ & $0$ & $0$
 & \hspace{-5pt}$-\p_{y_1}$\hspace{-5pt} & $-\p_{x_1}$ \\
 \hhline{|-||----|}
\vph $\p_{y_1}$ & $0$ & $0$ & \hspace{-5pt}$-\p_{x_1}$\hspace{-5pt}
 & $\p_{y_1}$ \\
 \hhline{|-||----|}
\vph $\p_{x_2}$ & $\p_{y_1}$  & $\p_{x_1}$  & $0$  &
 \hspace{-10pt}$-x_1\p_{y_1}$\hspace{-10pt} \\
 \hhline{|-||----|}
\vph $\p_{y_2}$ & \hspace{-5pt}$\p_{x_1}$\hspace{-5pt} &
 \hspace{-5pt}$-\p_{y_1}$\hspace{-5pt}
 & \hspace{-15pt} $x_1\p_{y_1}$\hspace{-15pt}  & $0$\\
 \hhline{-||----}
  \end{tabular}}
\hspace{8pt}
{\renewcommand{\arraystretch}{0}%
  \begin{tabular}{|c||c|c|c|c|}
 \hhline{-||----}
\vph
 \hspace{-5pt}$\nabla_{\!\!\scriptscriptstyle\uparrow}
 {\scriptstyle\leftarrow}$\hspace{-5pt}
 & $\p_{x_1}$ & $\p_{y_1}$ & $\p_{x_2}$ & $\p_{y_2}$ \\
 \hhline{:=::====:}
\vph $\p_{x_1}$ & $0$ & $0$ & $-\frac14\p_{y_1}$ & $-\frac34\p_{x_1}$ \\
 \hhline{|-||----|}
\vph $\p_{y_1}$ & $0$ & $0$ & $-\frac14\p_{x_1}$ & $-\frac14\p_{y_1}$ \\
 \hhline{|-||----|}
\vph $\p_{x_2}$ & $0$  & $0$  & $0$ & $0$ \\
 \hhline{|-||----|}
\vph $\p_{y_2}$ &
 \hspace{-4pt}$-\frac12\p_{x_1}$\hspace{-4pt} &
 \hspace{-4pt}$-\frac12\p_{y_1}$\hspace{-4pt} &
 \hspace{-12pt}$\frac14x_1\p_{y_1}$\hspace{-12pt} &
 \hspace{-12pt}$\frac14x_1\p_{x_1}$\hspace{-12pt} \\
 \hhline{-||----}
  \end{tabular}}
  \end{center}

So we find that the horizontal planes are
$H=\langle\p_{x_1},\p_{y_1}\rangle$, whence the structure on $N_\CC M$
is
 $$
\hat J\p_{x_1}=\p_{y_1},\
\hat J\p_{y_1}=-\p_{x_1},\
\hat J\p_{x_2}=\p_{y_2},\
\hat J\p_{y_2}=-\p_{x_2}
 $$
and $N_{\hat J}=0$.
 \end{ex}

 \begin{ex}
\underline{"Transversal distribution $\Pi^2$"}.
Let the structure be now
 $$
J\p_{x_1}=\p_{y_1}+x_2\p_{x_2},\
J\p_{y_1}=-\p_{x_1}-x_2\p_{y_2},\
J\p_{x_2}=\p_{y_2},\
J\p_{y_2}=-\p_{x_2}.
 $$
Again the curve $\CC=\{x_2=y_2=0\}$ is pseudoholomorphic and
the Nijenhuis tensor and a minimal almost complex connection are
(now $\Pi^2=\langle\p_{x_2},\p_{y_2}\rangle$):

 \vspace{1pt}
  \begin{center}
{\renewcommand{\arraystretch}{0}%
  \begin{tabular}{|c||c|c|c|c|}
 \hhline{-||----}
\vph
 \hspace{-5pt}$N_J({\scriptstyle\uparrow},
 {\scriptstyle\leftarrow})$\hspace{-5pt}
 & $\p_{x_1}$ & $\p_{y_1}$ & $\p_{x_2}$ & $\p_{y_2}$ \\
 \hhline{:=::====:}
\vph $\p_{x_1}$ & $0$ & \hspace{-10pt}$-x_2\p_{y_2}$\hspace{-10pt}
 & \hspace{-5pt}$\p_{y_2}$\hspace{-5pt} & $\p_{x_2}$ \\
 \hhline{|-||----|}
\vph $\p_{y_1}$ & \hspace{-15pt}$x_2\p_{y_2}$\hspace{-15pt} & $0$ &
 \hspace{-5pt}$\p_{x_2}$\hspace{-5pt} &
 \hspace{-5pt}$-\p_{y_2}$\hspace{-5pt} \\
 \hhline{|-||----|}
\vph $\p_{x_2}$ & \hspace{-5pt}$-\p_{y_2}$\hspace{-5pt}  & $-\p_{x_2}$
 & $0$  &  $0$ \\
 \hhline{|-||----|}
\vph $\p_{y_2}$ & \hspace{-5pt}$-\p_{x_2}$\hspace{-5pt} &
 \hspace{-5pt}$\p_{y_2}$\hspace{-5pt} & $0$ & $0$ \\
 \hhline{-||----}
  \end{tabular}}
\hspace{8pt}
{\renewcommand{\arraystretch}{0}%
  \begin{tabular}{|c||c|c|c|c|}
 \hhline{-||----}
\vph
 \hspace{-5pt}$\nabla_{\!\!\scriptscriptstyle\uparrow}
 {\scriptstyle\leftarrow}$\hspace{-5pt}
 & $\p_{x_1}$ & $\p_{y_1}$ & $\p_{x_2}$ & $\p_{y_2}$ \\
 \hhline{:=::====:}
\vph $\p_{x_1}$ & $0$ & $0$ & $0$ & $0$ \\
 \hhline{|-||----|}
\vph $\p_{y_1}$ &
 \hspace{-12pt}$\frac14x_2\p_{y_2}$\hspace{-12pt} &
 \hspace{-12pt}$\frac14x_2\p_{x_2}$\hspace{-12pt} &
 \hspace{-4pt}$-\frac12\p_{x_2}$\hspace{-4pt} &
 \hspace{-4pt}$-\frac12\p_{y_2}$\hspace{-4pt} \\
 \hhline{|-||----|}
\vph $\p_{x_2}$ & $-\frac14\p_{y_2}$ & $-\frac34\p_{x_2}$ & $0$ & $0$ \\
 \hhline{|-||----|}
\vph $\p_{y_2}$ & $-\frac14\p_{x_2}$ & $-\frac14\p_{y_2}$ & $0$ & $0$ \\
 \hhline{-||----}
  \end{tabular}}
  \end{center}

The horizontal planes are
$H=\langle\p_{x_1},
\p_{y_1}+\frac12e^{\frac12y_1}(x_2\p_{x_2}+y_2\p_{y_2})\rangle$.
So the structure on $N_\CC M$ is
 $$
\begin{array}{ll}
\hat J\p_{x_1}=\p_{y_1}+\frac12e^{\frac12y_1}(x_2\p_{x_2}+y_2\p_{y_2}), &
\hat J\p_{x_2}=\p_{y_2},\\
\hat J\p_{y_1}=-\p_{x_1}-\frac12e^{\frac12y_1}(x_2\p_{y_2}-y_2\p_{x_2}), &
\hat J\p_{y_2}=-\p_{x_2}.
\end{array}
 $$
Now $N_{\hat J}(\p_{x_1},\p_{x_2})=-e^{\frac12y_1}\p_{y_2}$ and so
the characteristic distribution of the normal structure is "the same"
as for the structure $J$:
${\hat\Pi}^2=\op{Im}N_{\hat J}=\langle\p_{x_2},\p_{y_2}\rangle$.

So the answer to the above question is positive if $N_J$ is of special
type.
 \end{ex}

\subsection{Linear bundle almost complex structures}

\hspace{13.5pt}
Consider an almost complex vector bundle
$\pi:(E,J)\stackrel{F}\to (\CC,J_0)$ of the rank $\dim F=2n$
and suppose the restriction $\left.J\right|_F$ is a linear complex
structure on the fiber. So we can also consider $(E,\pi,\CC)$ just as
vector bundle with complex structures in the fibers.

 \begin{dfn}\po
We call the almost complex structure $J$ on $E$
{\em linear bundle structure\/} if there exists a linear minimal
almost complex connection $\hat\nabla$ on this bundle
such that the lift $T_bE\stackrel{\ \hat\nabla}\leftarrow T_a\CC$
is a complex mapping, splitting the exact sequence
 $$
0\to F\to T_aE\to T_x\CC\to 0, \qquad x=\pi(a).
 $$
In particular the zero section $\CC\subset E$ is a $J$-pseudoholomorphic
curve.
 \end{dfn}

We note now that the almost complex structure $\hat J$ on the
normal bundle $N_\CC M$ is a linear bundle structure.

 \begin{lem}\po\label{LBACS}
The Nijenhuis tensor of linear bundle almost complex structure $J$
is constant along the fibers and determines $J$-invariant
differential system $\Pi=\op{Im}N_J\subset TE$ which is a subsystem
of vertical distribution $F$.
 \end{lem}

 \begin{proof}
There are local coordinates $(\vp,z)$ on $\pi^{-1}(U)=U\t F$,
$\vp=\vp_1+i\vp_2$, $z_k=x_k+iy_k$, $1\le k\le n$, such that
$\hat\nabla$-lift of $\p_{\vp_i}$ is
$\p_{\vp_i}+\sum b_{ij}\p_{x_j}+c_{ij}\p_{y_j}$, where
the coefficients are linear functions of $z_k$ and the vertical
coordinates are complex linear coordinates but horizontal coordinates are
complex only on the zero section. So the structure $J$ is given by
relations
 $$
J\p_{\vp_1}=\p_{\vp_2}+\sum
(b_{2j}+c_{1j})\p_{x_j}+(c_{2j}-b_{1j})\p_{y_j},
\qquad J\p_{x_k}=\p_{y_k}
 $$
and so $N_J(\p_{\vp_1},\p_{x_j})=\sum \a^k_j\p_{x_k}+\b^k_j\p_{y_k}$
with constant coefficients $\a^k_j,\b^k_j$.

The last claim follows because $N_J(F,F)=0$ and $\Pi=\C N_J(\p_{\vp_1},F)$
is a linear subspace of $F$ (of course invariant under $J$).
Note that $\op{rk}\Pi$ can vary with $\vp\in\CC$.
 \end{proof}

 \begin{rk}\po\label{rkonlft}
The conditions $N_J(F,F)=0$ and $\op{Im}N_J\subset F$
allow to lift naturally this tensor from $\CC$ to $E$.
Actually let $\x\in T_xC$ and $\e,\te\in T_xF\subset T_xE$.
Let us choose some lift $\ti\x\in T_aE$ (i.e.\ $(\pi_*)\ti\x=\x$).
Let $\ti{\vphantom{\x}\e},\ti\te\in T_aF=F$ be identified with $\e,\te$.
Then we define $N_J(\ti\x,\ti{\vphantom{\x}\e})=\ti\te$ at $a\in E$ if
$N_J(\x,\e)=\te$ at $x=\pi(a)$. Our assumptions imply that
this extension does not depend on a lift of $\x$.
So $N_J$ on $E$ is the canonical tensor.
 \end{rk}

Consider now the case $\op{dim}E=4$.
For every point $a\in E$ denote by $r=r_a\in T_aE$ the vertical
vector equal to $\vec{xa}\in F\simeq T_aF$ with $x=\pi(a)\in\CC$.

 \begin{theorem}\po\label{th:1}
Let $J$ be a linear bundle structure on a vector bundle $\pi:E\to\CC$
with 2-dimensional fibers over a curve $(\CC,J')$.
Then for some complex structure $J_0$ on $E$,
making the bundle holomorphic, we have:
 $$
J=J_0+\frac12J_0N_J(r,\cdot).
 $$
 \end{theorem}

 \begin{proof}
Let us define the structure by the formula
 \begin{equation}
J_0=J-\frac12JN_J(r,\cdot).
 \label{e:J}
 \end{equation}
Since $\left.N_J\right|_F\equiv0$ this structure $J_0|_F=J|_F$ is linear
on fibers. This proves the formula for $J$ of the theorem.

We first show that the structure $J_0$ is almost complex.
Note that by lemma~\ref{LBACS} $N_J(r,\x)\in F$ for any $\x$ and
$N_J(r,\x)=0$ for $\x\in F$. So
 $$
J_0^2=J^2-
\frac12J^2N_J(r,\cdot)-\frac12JN_J(r,J\cdot)+
\frac14JN_J(r,JN_J(r,\cdot))=J^2=-\1.
 $$

Now we show that this $J_0$ is integrable. By Newlander-Nirenberg
theorem \cite{NW} this is equivalent to the vanishing of the tensor
$N_{J_0}$. Let us choose local coordinates $(z,\vp)$ as in
proposition~\ref{prop:1}. In these coordinates $J\p_x=J_0\p_x=\p_y$ and we
deduce from~(\ref{e:J}):
 $$
 \begin{array}{rcl}
N_{J_0}(\p_x,\p_{\vp_1})
\!\!\!&=&\!\!\!
N_J(\p_x,\p_{\vp_1})-[\p_y,\frac12JN_J(x\p_x+y\p_y,\p_{\vp_1})]\\
\!\!\!&&\!\!\!
\hphantom{N\,(\p_x,\p_{\vp})}
\vphantom{\left\{\dfrac{1^1}{2^2}\right\}}
+J[\p_x,\frac12JN_J(x\p_x+y\p_y,\p_{\vp_1})]\\
\!\!\!&=&\!\!\!
N_J(\p_x,\p_{\vp_1})-\frac12N_J(\p_x,\p_{\vp_1})
-\frac12N_J(\p_x,\p_{\vp_1})=0.
 \end{array}
 $$
Since the bivector $\p_x\wedge\p_{\vp_1}$ generates $\Lambda^2_\C TE$
the claim follows.
 \end{proof}

Note that the tensor $N_J$ on the right-hand side in the formula of the
theorem can be defined only at the points of $\CC\subset E$ and then
lifted to $E$ as in remark~\ref{rkonlft}. Since one easily describes
complex structures on vector bundles (for trivial bundles see
\S\ref{HolRec}) we get the complete description of linear bundle almost
complex structures.

\subsection{Normal form of 1-jet of an almost complex structure}

Consider the ideal of real-valued functions corresponding to the curve
$\CC$
 $$
\mu=\{f\in C^\infty(M^4)\,|\,f(\CC)=0\}.
 $$
Degrees of this ideal give the filtration $\mu^k$ on every
$C^\infty(M)$-module, in particular we can talk about jets of
tensor fields.

 \begin{theorem}\po\label{th:2}
Let $\CC\subset M^4$ be a pseudoholomorphic curve
with respect to two almost complex structures $J_1$ and $J_2$ on $M^4$.
Assume $H^1(\CC;N_\CC M)=0$, i.e.\ the curve is a disk $D_R$, the entire
line $\C$ or the sphere $S^2$ with trivial normal bundle $S^2\cdot S^2=0$.
If $\left.J_1\right|_a=\left.J_2\right|_a$ and
$\left.N_{J_1}\right|_a=\left.N_{J_2}\right|_a$ at every point $a\in\CC$
then the structures $J_1$ and $J_2$ are 1-equivalent, i.e.\ there exists a
diffeomorphism $\psi$ of a neighborhood $\ok(\CC)$, preserving $\CC$,
such that
 $$
J_2=\psi^*J_1\,\op{mod}\mu^2.
 $$
 \end{theorem}

 \begin{proof}
This statement is an analog of theorem~1 \cite{K1}. Details of the
construction can be found there. Here we present a short proof with the
indication of differences.

By the hypothesis the desired diffeomorphism $\psi$ has 1-symbol
$\Phi^{(1)}=\op{id}\in\tau^*\ot\tau$ along the curve $\CC$,
where $\tau=TM|_\CC$. Its 2-symbol
$\Phi^{(2)}\in S^2\tau^*\ot\tau$ in local coordinates is given by
$\Phi^{(2)}(\x,\e)^i=\dfrac{\p^2\psi^i}{\p x^r\p x^s}\x^r\e^s$.
Moreover these symbols are compatible with the condition
$\psi|_\CC=\op{id}$:
 \begin{equation}
\x,\e\in T\CC\Rightarrow \Phi^{(2)}(\x,\e)=0.
 \label{e:COM}
 \end{equation}

The symbol satisfies the following condition:
 $$
dJ_1(\x,\e)+J_1\c\P^{(2)}(\x,\e)= \P^{(2)}(J_2\x,\e)+dJ_2(\x,\e),
 $$
where $d=d^\nabla$ is the differential w.r.t.\ some symmetric
(and so not almost complex) curvature-free connection $\nabla$.

Equivalently this means that the tensor
 \begin{equation}\label{e:P}
P(\x,\e)=J_1\c\P^{(2)}(\x,\e)-\P^{(2)}(J_2\x,\e)
 \end{equation}
satisfies the equation
 \begin{equation}\label{e:P'}
P(\x,\e)=dJ_2(\x,\e)-dJ_1(\x,\e).
 \end{equation}

Now given $J_1$ and $J_2$ we have $P$ and would like to find the
corresponding $\P$ from (\ref{e:P}). Form equality~(\ref{e:P'}) we deduce
 $$
J_1\c P(\x,\e)=-P(J_2\x,\e),
 $$
which implies $P(\x,\e)=J_1B(\x,\e)-B(J_2\x,\e)$ for some $(2,1)$-tensor
$B$. Now the identity
 \begin{equation}
P(\x,\e)-P(\e,\x)=P(J_2\x,J_2\e)-P(J_2\e,J_2\x)
 \label{e:ID}
 \end{equation}
implies the possibility to satisfy equation~(\ref{e:P}) by the choice
 $$
\Phi^{(2)}(\x,\e)=\frac12[B(\x,\e)+B(\e,\x)]
-\dfrac J4[B(\x,J\e)+B(\e,J\x)-B(J\x,\e)-B(J\e,\x)].
 $$

Identity~(\ref{e:ID}) means exactly $N_{J_1}=N_{J_2}$. Note that the
tensor $B$ can be chosen to satisfy~(\ref{e:COM}) and all the
constructions respect this condition ($T\CC$ is $J_k$-invariant).
So we get $\Phi^{(2)}$ satisfying~(\ref{e:COM}) and (\ref{e:P}).

Now the symbols $\Phi$ stand for the differential of the mapping $\psi$ we
seek for. Calculations above show the $\tau^*\ot\tau$-valued 1-form
generated by $\Phi$ is closed. Then the condition $H^1(\CC;N_\CC M)=0$
imply that it is exact. Hence the symbols $\Phi$ are
integrated to the diffeomorphism $\psi$ we sought for.
 \end{proof}

Consider PH-sphere $\CC=S^2\subset M^4$, $\CC\cdot\CC=0$.
On the family of transversal disks $D_\vp$, constructed in
proposition~\ref{prop:1}, we can choose smooth global
complex coordinate $z$ so that $z$ is complex w.r.t.\
the restriction $\left.J\right|_{D_\vp}$.
So a neighborhood $\ok$ of $S^2$ is represented as a product
$D^2\times S^2$ with coordinates $(z,\vp)$,
$\vp$ taking values in $\bar\C$.
Let us define a complex structure by the
product-formula $J_0=J_v\t J_h$, with the vertical part $J_v$ induced
by the projection along the spheres $S^2=\{z=\op{const}\}$
from the structure $J|_{D_\vp}$ and the horizontal part $J_h$ induced
by the projection along the disks $D_\vp$ from the structure $J|_\CC$.

 \begin{cor}\po\label{cor:131000}
Almost complex structure $J$ in a neighborhood $\ok\subset M^4$ of
PH-sphere $\CC=S^2$ is 1-equivalent to the almost complex structure
$J'=J_0+\frac12J_0N_J(r,\cdot)+\dots$, where $N_J$ is the field of the
Nijenhuis tensors along $\CC$. In other words there exists a local
diffeomorphism $\psi$, $\psi|_\CC=\op{id}$, such that
 $$
\psi^*J\x=J_0\x+\frac12J_0N_J(r_a,\x)\,\op{mod}\mu^2
 $$
for any vector $\x\in T_aM$, $a\in\ok$.
 \end{cor}

 \begin{proof}
The values of the almost complex structure $J$ and the almost complex
$\op{mod}\mu^2$ structure $\tilde J=J_0+\frac12J_0N_J(r_a,\cdot)+\dots$
on the right-hand side of the formula are equal at the points of the
sphere $S^2$. Let us show that the same is true for their Nijenhuis
tensors. Namely we show that the Nijenhuis tensor $N_{\tilde J}$
calculated according to definition (\ref{e0}) is equal
to the prescribed Nijenhuis tensor $N_J$ at all points $S^2$.

The calculations are local and it is sufficient to consider the value
of $N_{\ti J}$ on two complex independent vectors, say on $\p_x$ and
$\p_{\vp_1}$, where $\vp=\vp_1+i\vp_2$ and $z=x+iy$ are local
coordinates defining $J_0$. Let us denote by
$\doteq$ the equality $\op{mod}\mu$. We have:
 $$
 \begin{array}{rcl}
N_{\ti J}(\p_x,\p_{\vp_1})
\!\!\!&\doteq&\!\!\!
[J_0\p_x,\ti J\p_{\vp_1}]-J_0[\p_{x},\ti J\p_{\vp_1}]
-J_0[J_0\p_{x},\p_{\vp_1}]-[\p_{x},\p_{\vp_1}]
\\
\!\!\!&\doteq&\!\!\!
[\p_y,\p_{\vp_2}+\frac12J_0N_J(x\p_{x}+y\p_y,\p_{\vp_1})]
\\
\!\!\!&&\!\!\!
-J_0
[\p_{x},\p_{\vp_2}+\frac12J_0N_J(x\p_{x}+y\p_y,\p_{\vp_1})]
\\
\!\!\!&\doteq&\!\!\!
\frac12J_0N_J(\p_y,\p_{\vp_1})+
\frac12N_J(\p_{x},\p_{\vp_1})
\doteq
N_J(\p_{x},\p_{\vp_1}).
 \end{array}
 $$
Thus $N_{\ti J}=N_J$ along $S^2$ and the claim is proved.
 \end{proof}

 \begin{rk}\po
Another way to construct a complex structure $J'_0$ in a
neighborhood of PH-sphere $\CC=S^2\subset M^4$
with trivial normal bundle ($\CC\cdot\CC=0$) is to
foliate this neighborhood by homologous PH-spheres (\cite{M2})
and to induce the vertical part $J_v$ from some
particular transversal disks $D_{\vp_0}$.
 \end{rk}

Consider now arbitrary PH-curve $\CC$ with trivial normal bundle
$\CC\cdot\CC=0$.
Similarly to the proposition~\ref{prop:1} we find coordinates
$(z,\vp)$ in a neighborhood of $\CC\subset M$, where
$\vp$ is multivalued and $\{\vp=\op{const}\}$ is a family of transversal
PH-disks $D_\vp$. $D_\vp$ can be equipped with vector space structure
and so for every $a\in D_\vp\subset M^4$ we set
$r\in T_aD_\vp$ be the radius vector $\vec{0z}$ attached at the point $a$.
Here we use the canonical identification $T_aF\simeq F$ for the vector
spaces.

Theorem~\ref{th:1} expresses the linear bundle almost complex structure
in terms of some complex structure and a Nijenhuis tensor.
When almost complex structure is arbitrary then
the conclusions of proposition~\ref{lem:3} and theorem~1 are wrong
(the formula of theorem 1 can give the operator with $J^2\ne-\1$).
However on the level of 1-jet the theorem remains true.

 \begin{theorem}\po
 \label{th:3}
Let $\CC\subset (M^4,J)$ be a PH-curve and
$N_J\in\Lambda^2\tau^*\ot\tau$ be the field of Nijenhuis tensors
of $J$ along $\CC$, where $\tau=TM|_\CC$. Assume the
Nijenhuis tensor characteristic distribution $\Pi^2$ is transversal
everywhere to $\CC$. Then for some complex structure $J_0$ in a
neighborhood of $\CC$ we have:
 \begin{equation}
J\x=J_0\x+\frac12J_0N_J(r,\x)\,\op{mod}\mu^2.
 \label{e:MU}
 \end{equation}
 \end{theorem}

 \begin{proof}
We can assume the disks $D_\vp$ have $\Pi^2$ as tangent planes
at the points of $\CC$.
We use the symbol $\doteq$ as in the proof of the corollary above.
First note that formula (\ref{e:MU}) implies $J\doteq J_0$ for vertical
vectors (belonging to $TD_\vp$). Thus we define
$J_0\doteq J-\frac12JN_J(r,\x)$.
Similarly to the proof of theorem~\ref{th:1} it is proved that
$J_0^2\doteq-\1$ and that $N_{J_0}\doteq0$. Thus $J_0$ is a complex
structure $\op{mod}\mu^2$ \cite{K1} and we can take instead of it
any {\em complex\/} structure $\ti J_0=J_0\,\mod\mu^2$.
 \end{proof}

In this theorem $J_0$ is some complex structure while in the
corollary~\ref{cor:131000} we had some definite complex structure.
In the case $C=T^2$ we can also assume
the complex structure $J_0$ is $J_{(\l)}$ from proposition~\ref{prop:2}.

Actually, if the complex structure $J_0$ in the torus neighborhood has
normal bundle determined by the normal nonresonant pair $(\l,\nu)$,
then we can find coordinates $(z,\vp)$ in a
neighborhood with the gluing rule~(\ref{e1}) in which $J_0$ is standard.
In these coordinates
 \begin{equation}
J_\a^\g=\d_{\a}^{\g+n}+\frac12N_{\a\b}^{\g+n}r^\b+\bar{\bar o}(|r|),
 \label{*1}
 \end{equation}
with $z=r^2+ir^4$ and indices $1,3$ being used for the real and
imaginary part of $\vp=r^1+ir^3$. Here $N_{\a\b}^\g$ are the components of
the Nijenhuis tensor. We assume every object changes its sign when one
index is changed by $4$, e.g.\ $N_{\a\b}^{\g+4}=-N_{\a\b}^\g$.  The
summation index $\b$ is assumed to take only values 2,4.

 \begin{rk}\po
Similar statement holds also for 1-jet of almost complex structure
at a point. It says that for some coordinate system $(x^i)$ near
the point the operator $J$ has components
 \begin{equation}
J_\a^\g=\d_{\a}^{\g+n}+\frac14N_{\a\b}^{\g+n}x^\b+\bar{\bar o}(|x|),
 \label{*2}
 \end{equation}
where $N^\g_{\a\b}$ are components of the Nijenhuis tensor and we assume
$2n$-twisted periodicity, i.e.\ $x^{i+2n}=-x^i$ etc.

This formula is obtained by studying the structural function (Weyl tensor)
of the corresponding geometric structure (\cite{K1},\cite{KL}).
The difference of the factors (1/2 instead of 1/4)
in the formulas is connected with the fact that we consider
different coordinate systems: in (\ref{*1}) we consider the radial
vector field $r$ which connects the point to its projection on the torus
and in (\ref{*2}) we use the usual radius-vector directing from a
fixed origin. Such formulas for $J\,\op{mod}\mu^2$ can be considered along
any PH-submanifolds of $(M^{2n},J)$ but generically there are only 0- and
2-dimensional PH-submanifolds.
 \end{rk}

\section{Curves in a neighborhood of a PH-torus}

\subsection{Locally foliating families of cylinders}

\hspace{13.5pt}
In this section we consider only pseudoholomorphic tori $\CC=T^2$.
Results of section~1 imply there are moduli of almost
complex structures on a germ of $\CC$
contrary to the complex case in which Arnold found a
normal form of a neighborhood of an elliptic curve $T^2\subset M^4$.
To get some analogs of the complex situation we consider
foliations of a neighborhood of $T^2$ by pseudoholomorphic curves.
Generically there are no compact curves in a neighborhood and
PH-tori appear discretely because index of the
corresponding linearized Cauchy-Riemann operator is zero (\cite{Ku2}). So
one should consider noncompact curves. We begin with PH-lines.
Let $\|\!\cdot\!\|$ be some fixed norm.

 \begin{prop}\po \label{lem:4}
There are two neighborhoods $\ok'\subset\ok$ of $T^2\subset M^4$,
a change of almost complex structure in $\ok\setminus\ok'$ and a
number $C>0$ with the following properties.
For every $R>0$ there exists a smooth
family of PH-disks $f_\a:D_R\to\ok$ with uniformly bounded norms
$\|(f_\a)_*(z)\|\le C$ and $\|(f_\a)_*(0)\|=1$. Moreover this family
fills some smaller neighborhood $\ok''\subset\ok'$ of $T^2$, i.e.\
 $$
\ok''\subset\cup_\a f_\a(D_R).
 $$
 \end{prop}

 \begin{proof}
Let us take the universal covering $\hat\ok\simeq \C\t D^2$ of $\ok$. The
torus is covered by the entire line $\C\to T^2$. Changing the structure
$J$ at infinity in $\hat\ok$ and near the boundary
to the integrable one we glue the manifold to the product $S^2\t S^2$
with the line $\C$ being glued to the first factor $S^2_1$.
Then the introduction of the taming symplectic product-structure
$\oo=\oo_1\oplus\oo_2$ yields a foliation of $S^2_1\t S^2_2$ by PH-spheres
$S^2$ in the homology class of the first factor if we additionally demand
$\oo_1(S^2_1)<\oo_2(S^2_2)$, i.e.\ the homology class $[S^2_1]$ of the
first sphere-factor is symplectically simple. Here we use the fact that
the dimension is 4: due to positivity of intersections \cite{M1} we
actually have a foliation (\cite{M2}).

This foliation of
$S^2\t S^2$ gives a family of big PH-disks on $\hat\ok$ parametrized by
the radius $\rho$ of disk in $\C$ out of which the almost complex
structure is changed. To get estimates we use
Brody reparametrization lemma as in \cite{KO}.
The filling property is given by the choice of pertrubation
of the structure in $\ok\setminus\ok'$: This is simple if we require
the boundary of $\ok$ to be $J$-pseudoconvex (\cite{G} 2.4.D).
 \end{proof}

We now consider filling by pseudoholomorphic cylinders
$\cyl_R=[-R;R]\t S^1$. They can be considered also as annuli
$\cyl_R\subset\C\setminus\{0\}$, but
the next construction differs from the previous because the cylinders
are now non-contractible in $\ok$ (Fig.1).

 \begin{prop}\po \label{prop:3}
In the statement of proposition \ref{lem:4} we can change disks $D_R$
to $\cyl_R$ and get for every $R>0$ a family of
PH-cylinders $f_\a:\cyl_R\to\ok$ with uniformly bounded norms
and normalization $\|(f_\a)_*(0)\|=1$. In addition this family is
filling:
 $$
\ok''\subset\cup_\a f_\a(\cyl_R).
 $$
 \end{prop}

 \begin{picture}(300,160)
\put(-118,160){\centerline{{\special{em:graph bk1.gif}}}}
\put(85,10){Figure 1. Filling by PH-cylinders}
 \end{picture}


 \begin{proof}
Actually take a covering of the neighborhood $\ok$ which corresponds to
one cycle of the torus. The torus is covered by the entire cylinder
$\cli\to T^2$. We can change the almost complex structure $J$ at
infinity so that it makes possible to "pinch" each end of the cylinder.
This means we perturb the structure $J$ so that it is
standard integrable outside some $\cyl_{R_2}\subset\cli$ and the support is
also a big cylinder $\cyl_{R_1}$. Then we glue the ends to the disks.
This operation gives us a sphere $S^2$ instead of the cylinder
$\cli=\R\t S^1$. We can also assume that neighborhoods of two cylinder
ends are pinched (Fig.2).

 \begin{picture}(300,260)
\put(-142,240){\centerline{{\special{em:graph bk2.gif}}}}
\put(90,10){Figure 2. Cutting and Gluing}
 \end{picture}

Thus we have a neighborhood $U$ of the sphere $S^2_0$. It is
foliated by PH-spheres close to $S^2_0$. Actually, we can change the
structure $J$ near the boundary of this neighborhood, glue and get the
manifold-product $\hat M=S^2\t S^2$. As before it is foliated by
PH-spheres. Thus $U$ is foliated by PH-spheres and in the preimage they
give a PH-foliation by cylinders.
 \end{proof}

 \begin{rk}\po
Consider the infinite cylinder $\cli=\cyl_\infty$.
One can prove a more general statement:
There exists a change of the structure $J$ in $\ok\setminus\ok'$ such
that $\ok'$ is filled with images of entire PH-cylinders
$f_\a:\cli\to\ok$:  $\ok'\subset\cup_\a f_\a(\cli)$.
The technical difficulty is however that we need reparametrization in the
sequence of finite PH-cylinders passing the given point
and the limiting curve can shift and not pass through the point
(I am indebted to V.Bangert for a discussion about it).
To achieve filling one needs to control the sequence. We will consider the
problem elsewhere.
 \end{rk}

The positivity of intersections theory (\cite{M1}) is applicable only to
closed submanifolds and it cannot guarantee that leaves of
the considered filling do not intersect. If there are no
intersections of the leaves with close values of parameters and there are
no self-intersections we call the family {\em locally foliating\/}.

Recall that PH-torus $T^2$ is parametrized by periodic coordinate
$z\in\C$, i.e.\ the map $f_0:\C\to T^2$ satisfies
$f_0(z)=f_0(z+2\pi)=f_0(z+\nu)$, where $2\pi,\nu$ are periods. We
parametrize $f_\a$ by the condition $f_\a(z)=f_\a(z+2\pi)$.

Fix some transversal PH-disk $D_{\vp_0}$ to the torus $T^2$.
Consider the return map $\Psi:D_{\vp_0}\to D_{\vp_0}$ for our family of
surfaces $f_\a$. This is defined as follows: if a point $z\in D_{\vp_0}$
belongs to $f_\a$, $z=f_\a(\vp)$ we set $\Psi(z)=z'$ with
$z'=f_\a(\vp+\z(\vp))$, $\z(\vp)\approx\nu$, chosen to belong to
$D_{\vp_0}$ again. This
return map $\Psi$ is determined non-uniquely since several cylinders of
the family $f_\a$ can pass through $z$.

 \begin{prop}\po
Let the return map be either attracting or repelling,
$|\Psi(z)/z|\in(0,\ve)\text{ or }(\ve^{-1},\infty)$ for all $z\in
D_{\vp_0}$ close to zero and some fixed $\ve<1$ (we assume the inequality
holds notwithstanding the nonuniqueness of $\Psi$). Then the constructed
family $f_\a:\cli\to\ok$ is locally foliating in a neighborhood of $T^2$.
 \label{pro:4}
 \end{prop}

 \begin{proof}
The condition that $\Psi$ is attracting or repelling means that there
are no self-intersections of the leaves $f_\a(\cli)$.
If a leaf $f_\a(\cli)$ transversally intersect another $f_\b(\cli)$
with $\a\approx\b$ then because of the convergence
$f_\a(\cyl_R)\to f_\a(\cli)$ there is an intersection of the spheres
$S^2\subset \hat M$ from which cylinders are constructed. This is
certainly impossible because our family of spheres is foliating.
If the leaves $f_\a(\cli)$ and $f_\b(\cli)$ are tangent let us consider the
first order of their jets which are different.  This is the next number
after the tangency order.  Since the maps of finite cylinders $\cyl_R$ tends
to the maps of $\cli$ in $C^\infty$-topology the spheres in the manifold
$\hat M$ do intersect.  Every such intersection even nontransversal
contributes positively to the intersection number (\cite{M1}) which
contradicts $[S^2]\cdot[S^2]=0$. It remains to note that if two
pseudoholomorphic curves in an almost complex manifold have the infinite
tangency they must coincide (\cite{MS}).
 \end{proof}

For arbitrary almost complex structures it is not clear if there are
foliations of the $T^2$-neighborhood by cylinders.
The following can be probably solved using a method
similar to the Moser's proof~\cite{Mo} of PH-lines foliation
persistence.

{\bf Question:} Let $T^2$ be normal nonresonant elliptic
curve.  Let $J$ be a small almost complex perturbation of the complex
structure $J_0$. The curve $T^2$ will be perturbed into close PH-curve
$\CC$.  Is it true that a neighborhood $\ok(\CC)$ is foliated by
PH-cylinders?

\subsection{Floquet theory for the PH-foliations}

\hspace{13.5pt}
Holomorphic bundle (\ref{e1}) possesses a distinguished foliation
$\{z=\op{const}\}$. We call this foliation the {\em foliation by
$\l$-twisted cylinders\/}. Arnold's "Floquet-type" result in \cite{A1}\S27
implies that a neighborhood of every elliptic curve has a foliation by
twisted cylinders which is biholomorphic to the standard one. Here we
study PH-foliations.

Let $f_\a:\BB\to\ok$ be a foliating family of a neighborhood of a
PH-curve $\CC$ with trivial self-intersection. Let
$D_\vp=\{\vp=\op{const}\}$ be a family of normal disks from
proposition~\ref{poravstavat'}.
Then every path $\g(t)$ on $\CC$ with $\g(0)=\vp_0$,
$\g(1)=\vp_1$ gives a mapping $\Phi_\g:D_{\vp_0}\to D_{\vp_1}$
of shift along the leaves of $f_\a$. For a loop
$\g$ we have an automorphism of $D_{\vp}$. Since $f_\a$ is foliation there
is no local holonomy: $\Phi_\g=\op{id}$ for contractible loops $\g$.
Thus we can consider the map $\pi_1(\CC)\to\op{Aut}(D_\vp)$.

 \begin{dfn}\po
We call $\Phi_\g\in\op{Aut}(D_\vp)$ the {\it monodromy map\/} along
$\g\in\pi_1(\CC)$.
 \end{dfn}

For example there is no monodromy for the sphere $\CC=S^2$ and each
choice of local coordinates in a normal disk $D_{\vp_0}$ gives coordinates
for the others $D_\vp$.

Let now $\CC=T^2(2\pi,\nu)$. Since our foliating family $f_\a$ consists of
cylinders there is no monodromy along one generating cycle, let along
the cycle $\vp\mapsto \vp+2\pi$. Denote by $\Phi_\nu$ the monodromy
along the other cycle $\vp\mapsto\vp+\nu$.

Let $f_\a:\cli\to\ok$ be a foliating family of PH-cylinders
in a neighborhood $\ok'\subset\ok$ of the PH-torus $T^2$.
Assume the monodromy of the family $f_\a$ is
 \begin{equation}
\Phi_\nu(z)=\l z+\bar{\bar o}(|z|).
 \label{e:la}
 \end{equation}

Recall~\cite{A1} that $\l\in\C\setminus\{0\}$ is of $(C,\z)$-type if
$|\l^m-1|\ge C/m^{1+\z}$ for all $m\in\Z_{>1}$ (in particular
such are all nonzero from Poincar\'e domain $|\l|\ne1$).

 \begin{theorem}\po\label{th:4}
Let the monodromy $\Phi_\nu$ be a germ of biholomorphic mapping with
the number $\l$ (\ref{e:la}) of $(C,\z)$-type.
Then the foliating family of cylinders $f_\a$
is similar to $\l$-twisted foliation in the following sense:
There exist coordinates $(z,\vp)$ with the gluing rule (\ref{e1})
such that transversal disks $D_\vp=\{\vp=\op{const}\}$
are pseudoholomorphic, $\vp$ is a complex coordinate on the
curve $\CC$, $z$ is a complex coordinate on one $D_{\vp_0}$ and
the foliation $f_\a$ is given by $\{z=\op{const}\}$.
 \end{theorem}

 \begin{proof}
A coordinate system $z$ on a transversal disk $D_{\vp_0}$
provides coordinates on all others $D_\vp$.
These coordinates are multivalued because rotation along the
second cycle $\vp\mapsto\vp+\nu$ yields the monodromy
mapping $z\mapsto\psi(z)=\l z+\bar{\bar o}(|z|)$.

Now we apply Poncar\'e-Siegel theorem (\cite{A1}) in complex dimension 1
which states that the monodromy map is conjugate to the map
$\psi(z)=\l z$.
 \end{proof}

\subsection{Monodromy and transports}\label{S3.3}

\hspace{13.5pt}
Unlike complex case almost complex monodromy can be
non-holomorphic mapping of the fibers. Actually there are simple
examples of PH-foliations with any prescribed monodromy.

Moreover even if the monodromy is complex
(as required in theorem~\ref{th:4})
the transport maps $\Phi_\g:(D_{\vp_0},J)\to(D_{\vp_1},J)$
can be not. In fact this is the occasion of $\op{codim}=\infty$.

 \begin{prop}\po\label{prop:5}
Let $\CC\subset M$ be a PH-curve in a 4-dimensional manifold and
let $f_\a:P\to\ok$ be a local PH-foliating family of some
neighborhood $\ok(\CC)$. Let also $D_\vp$ be a transversal
PH-foliation and $\Phi_\g$ be the corresponding transport map
$D_{\vp_0}\to D_{\vp_1}$ along curves $\g(t)\subset\CC$,
$\g(0)=\vp_0$, $\g(1)=\vp_1$.
If all $\Phi_\g$ are holomorphic
then the Nijenhuis tensor characteristic distribution $\Pi^2$
is tangent to the leaves of $f_\a$.
 \end{prop}

 \begin{proof}
Actually the foliation provides a local bundle
$\pi:\ok\to D_\vp$. The hypothesys that all transports are complex
is equivalent to the claim that the bundle $\pi$ with fibers $f_\a(P)$
is almost complex. Therefore the statement follows from
proposition~\ref{lem:3}.
 \end{proof}

 \begin{rk}\po
If the distribution $\Pi^2$ is not integrable the transports are
not complex. But even integrability is not sufficient, see \S\ref{S2.2}.
 \end{rk}

\subsection{PH-tori deformation problem}\label{3.4}

\hspace{13.5pt}
Generically holomorphic 2-torus in a complex manifold cannot be deformed
(\cite{A1}). The same situation is also with PH-tori in almost complex
manifolds. This follows from vanishing of the index of the linearized
Cauchy-Riemann operator (\cite{Ku1}). By the deformation we mean
existence of close homologous PH-torus of the same periods
$T^2=T^2(2\pi,\nu)$. In this section we consider some examples where we
can make the condition of non-existence explicit.

\bff 1
Let us consider linear bundle almost complex structure $J$ on the bundle
$E\to T^2$. There exist coordinates $(z,\vp)$ with the gluing
rule (\ref{e1}) such that
 $$
\left\{
 \begin{array}{ll}
J\p_x=\p_y, & J\p_{\vp_1}=\p_{\vp_2}+x\cdot v-y\cdot Jv,
\\
J\p_y=-\p_x, & J\p_{\vp_2}=-\p_{\vp_1}-x\cdot Jv-y\cdot v.
 \end{array}
\right.
 $$
This formula follows from theorem~\ref{th:1} and the coordinates are
determined by $J_0$. Vector field $v=\frac12JN_J(\p_x,\p_{\vp_1})$ can be
decomposed $v=\a\p_x+\b\p_y$ with $\a=\a(\vp)$, $\b=\b(\vp)$. The
complexified vector bundle is decomposed $T_{\footnotesize\C}E=E_++E_-$,
where $E_\pm=\{\x\,|\,J\x=\pm i\x\}$; $E_-=\bar E_+$. Vectors
 $$
U_1=\p_{\vp}-z\bar b\,\p_{\bar z},\quad
U_2=\p_z,
 $$
form a basis of $E_+$. Here $\p_{\vp}=\frac12(\p_{\vp_1}-i\p_{\vp_2})$,
$\p_z=\frac12(\p_x-i\p_y)$ and $\bar b=\dfrac{\b+i\a}2$.
Thus the basis of $E_+^*$ in the decomposition
$T_{\footnotesize\C}^*E=E_+^*+E_-^*$ ($E_-^*=\bar E_+^*$) is
 $$
\oo_1=d\vp,\quad
\oo_2=dz+\bar zb\,d\bar\vp.
 $$

Now every pseudoholomorphic torus in $E$ homologous to the zero section
$T^2$ is of the form $f(T^2)$ for some section $f$ of the bundle
$E\to T^2$. This is to say each PH-torus in $E$ has unique transversal
intersection with every fiber. Actually we may compactify the fibers of
the bundle to the spheres and the claim follows from the positivity of
intersections (or even simpler by studying the degree of the projection
of this torus to the torus-base).

Let us deduce the equation for $f$. The curve $f(T^2)$ is pseudoholomorphic
iff
 $$
\left.\oo_2\right|_{z=f(\vp)}=c\cdot\oo_1.
 $$
Substituting $df=f_\vp d\vp+f_{\bar \vp}d\bar \vp$ we have:
 \begin{equation}
f_{\bar\vp}+ b\bar f=0.
\label{e:f}
 \end{equation}

 \begin{theorem}\po\label{chi}
Let $J$ be a linear bundle almost complex structure and $J_0$ be the
corresponding complex structure from the decomposition of
theorem~\ref{th:1}. Suppose the number $\l$, determined by the
complex structure $J_0$ in the bundle $E$ via {\rm (\ref{e1})},
is of unit length: $|\l|=1$.
Assume also that the function $\Lambda\in C^\infty(T^2)$,
determined uniquely by the equation
$\frac12JN_J(\p_z,\p_\vp)=\Lambda\p_{\bar z}$,
is nonzero holomorphic function: $\p_{\bar\vp}\Lambda=0$, $\Lambda\ne0$.
Then zero section $T^2$ is the unique PH-torus in $E$.
 \label{th:last}
 \end{theorem}

 \begin{proof}
First note that since
  \begin{equation}\label{beom}
 \begin{array}{ll}
\ddfrac12JN_J(\p_z,\p_\vp)=(\a-i\b)\p_{\bar z},
&\ddfrac12JN_J(\p_{\bar z},\p_\vp)=0,\\
\ddfrac12JN_J(\p_z,\p_{\bar\vp})=0,
&\ddfrac12JN_J(\p_{\bar z},\p_{\bar\vp})=(\a+i\b)\p_z,
 \end{array}
  \end{equation}
we have $\Lambda=-2i\bar b$. Thus
$\Lambda_{\bar\vp}=0\Leftrightarrow{\bar\Lambda}_\vp=0
\Leftrightarrow b_\vp=0$.

Let us show the equation (\ref{e:f}) has no nonzero solutions.
Complex Laplacian of $f$ equals
$f_{\bar\vp\vp}=-b\bar f_\vp=-b\overline{f_{\bar\vp}}=|b|^2f$.
Our torus neighborhood is the trivial cylinder
${\cal C}^2=\{\vp\in\C\,|\,\op{Im}\vp\in[0,\op{Im}\nu)\}/2\pi\Z$
neighborhood glued by the rule
$(z,\vp)\mapsto(\l z,\vp+\nu)$. Thus when $\vp\mapsto\vp+\nu$ we have:
$f\mapsto\l f$. So integrating over the cylinder gives:
 \begin{eqnarray}
\int_{{\cal C}^2}(f_{\bar\vp\vp}\bar f+f_{\bar\vp}\overline{f_{\bar\vp}})
\dfrac i2d\vp\we d\bar\vp
&\!\!\!=\!\!\!&
\int_{{\cal C}^2}\dfrac\p{\p\vp}(f_{\bar\vp}\bar f)
\dfrac i2d\vp\we d\bar\vp=
\int_{{\cal C}^2}\dfrac i2 d(f_{\bar\vp}\bar f d\bar\vp)=
\nonumber\\
&\!\!\!=\!\!\!&
(\l\bar\l-1)\oint_{S^1}\dfrac i2 f_{\bar\vp}\bar f d\bar\vp=0.
\nonumber
 \end{eqnarray}
So using the calculation with the Laplacian we have:
 \begin{equation}
\int_{{\cal C}^2}(|b|^2|f|^2+|f_{\bar\vp}|^2)d\vp_1\we d\vp_2=0.
 \label{e:int}
 \end{equation}
Therefore since $|b|\ne0$ we get $f=0$. Thus there are no homologous to
the zero section PH-tori $\ti T^2$ with $f\ne0$. If the homology class
of $\ti T^2$ is a multiple of the zero section $[\ti T^2]=k[T^2]$ a
$k$-finite covering finishes the proof.
 \end{proof}

 \begin{rk}\po
If $|b|=0$, i.e.\ almost complex structure $J$ is integrable $J=J_0$,
equality (\ref{e:int}) implies that $f$ is holomorphic section. Thus if
$\l^n\ne1$ we get again $f=0$ comparing the Fourier coefficients of $f$.
 \end{rk}

\bff 2
Consider a general almost complex structure $J$ with Nijenhuis tensor
characteristic distribution $\Pi^2$ transversal to some
PH-torus $T^2$. The linearized equation for close PH-tori
can be written in the form
 \begin{equation}\label{e:ff}
f_{\bar\vp}+af+b\bar f=0.
 \end{equation}
Actually, the linearization does not depend on a change of the structure
$J$ by second order quantities. Thus we can perturb $J$ to make the
distribution $\Pi^2$ integrable in $\ok\supset T^2$. This new almost
complex structure is given by the formula (\ref{251000}).

Let's write the equation for close PH-tori. The basis of $E_+$ is
 $$
U_1=\p_\vp+\dfrac{\bar A}{A+2i}\p_{\bar\vp}+\dfrac{\bar B}{A+2i}\p_{\bar
z},\quad U_2=\p_z,
 $$
where $A=A_1+iA_2$, $B=B_1+iB_2$. The corresponding basis of
$E_+^*$ is
 $$
\oo_1=d\vp-\dfrac{A}{\vphantom{A^{A'}}\bar A-2i}d\bar\vp,\quad
\oo_2=dz+\dfrac14\dfrac{\bar A B}{A_2+1}d\vp-
\dfrac14\dfrac{(A+2i)B}{A_2+1}d\bar\vp.
 $$
So the equation $\left.(\oo_2-c\cdot\oo_1)\right|_{z=f(\vp)}=0$
implies the required equation on $f$:
 \begin{equation}\label{eq_for_f}
f_{\bar\vp}+\dfrac{A}{\vphantom{A^{A'}}\bar A-2i}f_\vp-
\dfrac{B}{\vphantom{A^{A'}}\bar A-2i}=0.
 \end{equation}
Denote by $A^0$ and $B^0$ linearizations by fiber coordinate of the
functions $A$ and $B$ respectively.
Note that linearization of equation~(\ref{e:ADD}) implies that
$A^0$ is holomorphic w.r.t.\ $z$, that is we can bring our equations
to the constant $A^0$.

Since $\left.A\right|_{z=0}=0$, $\left.B\right|_{z=0}=0$,
linearization of equation (\ref{eq_for_f}) is
 $$
f_{\bar\vp}-\dfrac i2 B^0(f)=0,
 $$
which has the form (\ref{e:ff}) if we set $-\frac i2 B^0=az+b\bar z$.

Since we have equations (\ref{e:ADD}) the Nijenhuis tensor of
(\ref{251000}) is
 \begin{align*}
N_J(\p_x,\p_{\vp_1})&=
\left(\dfrac{\p B_1}{\p y}-\dfrac{\p A_1}{\p x}B_1+
\dfrac{\p A_2}{\p x}\dfrac{A_1B_1-B_2}{A_2+1}+\dfrac{\p B_2}{\p x}
\right)\p_x\\
&+
\left(\dfrac{\p B_2}{\p y}-\dfrac{\p A_1}{\p x}B_2+
\dfrac{\p A_2}{\p x}\dfrac{A_1B_2+B_1}{A_2+1}-\dfrac{\p B_1}{\p x}
\right)\p_y.
 \end{align*}
Therefore linearizing $A$ and $B$ we conclude that its values on $T^2$ are
given by the formula (see also (\ref{beom}))
 $$
\frac12JN_J(\p_z,\p_\vp)=-2i\bar b\,\p_{\bar z}.
 $$

 \begin{rk}\po
Since the only invariant of 1-jet of $J$ on a PH-curve is the Nijenhuis
tensor, which we expressed by the function $b(\vp)$, we can bring the
function $a(\vp)$ in (\ref{e:ff}) to the simplest form. Namely we can
introduce coordinates $(\vp,z)$ using the complex structure $J_0$ of
theorem~\ref{th:3}. This gives $a=0$ for normal coordinate $z$ on the
torus with gluing rule (\ref{e1}). Alternatively we can have global
well-defined coordinate $z$ but then $a=\op{const}$. This
proves a suggestion on p.\,430 \cite{Mo} that {\em "the linearized
equation can be brought into the form (\ref{e:ff}) with $a=\op{const}$"\/}.
 \end{rk}

 \begin{theorem}\po
Let almost complex structure $J$ in a neighborhood of PH-curve $T^2$ be
described by formula (\ref{e:MU}) with complex structure $J_0$ having
$|\l|=1$.  If the characteristic distribution $\Pi^2$ is transversal to
$T^2$ and for linearized structure $b(\vp)$ is anti-holomorphic,
then the curve $T^2$ is isolated and persistent.
 \end{theorem}

 \begin{proof}
Actually as Moser~\cite{Mo} noticed if the linearized equation
$f_{\bar\vp}+af+b\bar f=g$ has a unique solution for any
$g\in C^\infty(T^2;\C)$ then the torus is isolated and persistent.
But the linearization we studied in theorem \ref{chi}.
 \end{proof}

\bff 3
Note that in holomorphic bundle with $\l^n=1$ one can find tori $f(T^2)$
with $f\ne0$ of the type $T^2(2\pi kn,\nu l)$, which cover zero section
torus $T^2=T^2(2\pi,\nu)$. But in a fixed homology class
all PH-tori are of the same holomorphic type:

 \begin{lem}\po
Let $\ti T^2_1,\ti T^2_2\subset E$ be two PH-tori in a linear almost
complex bundle $E\to T^2$. If they are homologous then they are
biholomorphic.
 \end{lem}

 \begin{proof}
First consider tori in the homology class of the zero section $T^2$. As
was shown before theorem~\ref{th:last} the projection is a diffeomorphism.
Since in linear almost complex bundles $E\to T^2$ the projection is an
almost complex mapping, its restriction is a biholomorphism of the tori:
$\ti T^2_1\simeq T^2\simeq\ti T^2_2$. The general case
$[\ti T_i^2]=k[T^2]$ follows from the case $k=1$ by means of a
$k$-covering.
 \end{proof}

If we do not demand the bundle condition the opposite situation can
occur: in example \cite{A1}\S27 a neighborhood of the torus is foliated by
holomorphic tori of different holomorphic type. Similar situation occurs
also in almost complex case and invariants of section 1 can be nontrivial:

 \begin{ex}
Consider a foliation $f_\a:T^2\to\ok$ of 4-dim neighborhood of some
torus. Introduce the structure $J$ in horizontal directions
so that all the tori $T^2_\a$ are pseudoholomorphic but nonequivalent
(the parameter $\nu$ is changing).
Choose the structure $J$ on the normals $D_\vp$ so that the
transports $\Phi_\g$ are not holomorphic (\S\ref{S3.3}).
Define $J$ globally by the product formula. Then the distribution
$\Pi^2$ is nonintegrable and we get the distribution
$L^1=\Pi^3\cap T(T^2)$ (possibly with singularities).
 \end{ex}

\bff 4
Note that the tori in the same homology class can occur both in families
and discretely.

 \begin{ex}
Let $(T^4,J_0)$ be the standard complex torus, i.e.\ quotient of $\C^2$ by
the lattice $\Z^4$. Consider a PH-torus $T^2_0\subset T^4$. It is possible
to perturb $J_0$ in a neighborhood $\ok$ of $T^2_0$ so that the new
structure $J$ is isomorphic to the model structure near the torus
$T^2(2\pi,\nu)$ given by (\ref{e1}) in a neighborhood $\ok'\subset\ok$ and
$J=J_0$ outside $\ok$. Then there is 2-parametric family of PH-tori of
$[T^2_0]$-homology class outside $\ok$ and a unique PH-torus $T^2_0$ inside
$\ok'$.
 \end{ex}

\bff 5
Note that for PH-torus $\CC=T^2\subset (M^4,J)$ with $N_J|_a=0$ for all
$a\in T^2$ the normal bundle $N_\CC M$ is holomorphic.

 \begin{prop}\po
If the Nijenhuis tensor vanishes along a PH-torus and the pair $(\l,\nu)$,
characterizing the holomorphic bundle $N_\CC M$, is nonresonant then
small neighborhood $\ok$ of this torus cannot be foliated by
$T^2(2\pi,\nu)$-tori.
 \end{prop}

 \begin{proof}
Actually if there is a PH-foliation by tori then the linearization of
this foliation determines a holomorphic foliation of the normal bundle
which is impossible by~\cite{A1}\S27.
 \end{proof}

\bff 6
An intermediate condition on foliation between holomorphic and
pseudoholomorphic is that it be pseudoholomorphic with complex
transports. Proposition~\ref{prop:5} implies

 \begin{prop}\po
If the distribution $\Pi^2$ in a neighborhood $\ok\supset T^2$ is not
integrable or is integrable with noncompact leaves, then $\ok$ cannot be
foliated by PH-tori with complex transports.
 \end{prop}

Note that in the considerations above we needed to fix a holomorphic
structures on the tori sought for. The last proposition does not require
this.

\appendix
\section{Normal bundle in Riemannian geometry.}

\hspace{13.5pt}
The construction of structure on the normal to the submanifold bundle from
\S\ref{AC on NB} can be carried out for some other cases.

For example in Riemannian geometry there exists a unique Riemannian
connection $\nabla$. This Levi-Civita connection splits the normal bundle
$N_LM$ of any submanifold $i:L\subset M$ with induced Riemannian metric
$g_L=i^*(g_M)$. However the uniqueness of the connection $\nabla$ makes
the situation more rigid and for our construction we should require that
$L$ is a totally geodesic submanifold.

So we construct the metric in $N_LM$ and then again ask about relations
between two tensors: Riemannian curvature $R_g$ of the manifold $M$ at
points of $L\subset M$ and the curvature ${\hat R}_g$ for the total space
of $N_LM$ at zero section $L\subset N_LM$. In general there are no
relations, but for some parts of the tensors there is.

Namely consider the curvature of the normal bundle $R^\perp$,
which is the curvature tensor of the normal connection
$\nabla^\perp$ given by the orthogonal decomposition in
$TM|_L=TL\oplus N_LM$, $W=W_{\|}+W_\perp$.
Note that $R^\perp(X,Y)={\hat R}_g(X,Y)$ for $X,Y\in TL$ and the left-hand
side is not defined for others $X,Y$.

Let ${\rm II}:TL\ot TL\to N_LM$ be the second quadratic form of $L$ and
$A:TL\ot N_LM\to TL$ be the shape (Peterson) operator given by
$g(A(X,V),Y)=g({\rm II}(X,Y),V)$, $X,Y\in TL$, $V\in N_LM$.
Then if we denote $R^L$ the curvature of the Levi-Civita connection of $L$
we have the following famous equations:
 $$
 \begin{array}{llr}
\hspace{-14pt}
&[R_g(X,Y)Z]_{\|}=R^L(X,Y)Z+A(Y,{\rm II}(X,Z))-A(X,{\rm II}(Y,Z))
&\text{(Gauss eq.),} \\
\hspace{-14pt}
&[R_g(X,Y)Z]_{\perp}=(\nabla_Y{\rm II})(X,Z)-(\nabla_X{\rm II})(Y,Z)
&\hspace{-50pt}\text{(Codazzi-Maynardi eq.),} \\
\hspace{-14pt}
&[R_g(X,Y)V]_{\perp}=R^\perp(X,Y)V+{\rm II}(X,A(Y,V))-{\rm II}(Y,A(X,V))
&\text{(Ricci eq.),}
 \end{array}
 $$
where $X,Y,Z\in TL$, $V\in N_LM$.

In particular when $L$ is totally geodesic ${\rm II}=0$ and $A=0$, so that
the equations above mean $R_g(X,Y)={\hat R}(X,Y)$ for $X,Y\in TL$ at the
points of $L$. However there are no relations if we allow
$X,Y$ to be arbitrary from $TM$.

So horizontal parts of the both curvatures $R_g$ and ${\hat R}_g$
coincide. Note that in the almost complex case this is trivially so
because these horizontal parts vanish.

{\bf Question:} What happens in other geometries -- conformal, projective
etc? There are also notions of normal connections (\cite{N})
but on the normal spaces considered as bundles not manifolds.

\vspace{10pt}

{\footnotesize Acknowledgements.
I would like to thank V.\,I.\ Arnold who addressed me with several
questions and suggestions originating the paper. During the work I had the
benefit of discussions Riemannian geometry with D.\,V.\ Alekseevskij and
projective geometry with V.\ Goldberg.
I also express my thanks to V.\ Bangert for inviting me to the University
of Freiburg and for valuable discussions during the visit.}

\ {\hbox to 10.5cm{ \hrulefill }}

{\footnotesize
Note on figures: if they are not visible look at
www.math.uit.no/seminar/preprints.html}
%
%
%

\ {\hbox to 10.5cm{ \hrulefill }}
\bigskip

{\footnotesize
\hspace{-19pt}
Department of Math.\ and Stat., University of Tromsoe, Norway

\hspace{-19pt}
Mathem.\ Modelling Chair, Moscow Baumann State
Technological University, Russia}

{\it \hspace{-19pt} E-mail:} \quad
{\footnotesize kruglikov\verb"@"math.uit.no}

\end{document}